\documentclass[12pt]{amsart}
\usepackage{amsmath,mathrsfs}
\usepackage{amssymb}
\usepackage{latexsym}
\usepackage{amsfonts}

\newcommand{\w}{\omega}
\newcommand{\la}{\lambda}
\newtheorem{Pa}{Paper}[section]
\newtheorem{Tm}[Pa]{{\bf Theorem}}
\newtheorem{La}[Pa]{{\bf Lemma}}

\newtheorem{Rk}[Pa]{{\bf Remark}}
\newtheorem{Pn}[Pa]{{\bf Proposition}}

\newtheorem{Dn}[Pa]{{\bf Definition}}

\begin{document}
\title{Linear stochastic systems:
A white noise approach}
\author[D. Alpay]{Daniel Alpay}
\address{(DA) Department of mathematics,
Ben-Gurion University of the Negev, P.O. Box
653, Beer-Sheva 84105, Israel}
\email{dany@math.bgu.ac.il}
\author[D. Levanony]{David Levanony}
\address{(DL) Department of electrical
engineering, Ben-Gurion University of the Negev,
P.O. Box 653, Beer-Sheva 84105, Israel}
\email{levanony@ee.bgu.ac.il}
%

\date{}

\begin{abstract}
Using the white noise setting, in particular the
Wick product, the Hermite transform, and the
Kondratiev space, we present a new approach to
study linear stochastic systems, where
randomness is also included in the transfer
function. We prove BIBO type stability theorems
for these systems, both in the discrete and
continuous time cases. We also consider the case
of dissipative systems for both discrete and
continuous time systems.
We further study $\ell_1$-$\ell_2$ stability in
the discrete time case, and ${\mathbf
L}_2$-${\mathbf L}_\infty$ stability in the
continuous time case.
\end{abstract}

\maketitle
\tableofcontents

\section{Introduction}
\setcounter{equation}{0}

In this paper we propose a new approach for the
study of uncertainty within the theory of linear
stochastic systems, and prove a number of
stability theorems. To set the problems and
results in perspective we begin with a brief
historical introduction. Linear system theory,
operator theory and the theory of analytic
functions have a long history of interactions,
and two notable milestones are the work of M. Livsic on
the characteristic operator function, see
\cite{MR9:446c}, \cite{livsic}, \cite{l1}, and
the work of R. Kalman, see \cite{MR27:2147}. The
discussion of what is linear system theory would
lead us too far away, and we refer the reader to
\cite{MR569473} for more information. We also
refer to \cite{heltonbook} and
\cite{MR2002b:47144} for surveys and to
\cite{MR1200235} for a discussion of stability
results in the continuous time case.\\

For the purpose of this introduction, a
discrete-time, time-invariant linear system will be described by
an input-output relation of the form
\begin{equation}
\label{patricia_kaas}
y_n=(h*u)_n=
\sum_{m\in{\mathbb Z}} h_{n-m}u_m,\quad
n\in{\mathbb Z}.
\end{equation}
In this expression, the $h_n$ are pre-assigned
complex numbers, which stand for the impulse
response of the underlying system, and the input
and output are required to define a continuous
map between specified spaces ${\mathcal H}_i$
and ${\mathcal H}_o$ (and of course, this entails
conditions on the coefficients $h_n$). These
various conditions are translated into
properties for the ${\mathcal Z}$-transform
\begin{equation}
\label{Ztr}
\widehat{h}(\zeta)=\sum_{\mathbb Z}\zeta^nh_n
\end{equation}
of the sequence $(h_n)$. For instance, when the
system \eqref{patricia_kaas} is causal, that is,
when $h_n=0$ for $n<0$, it defines a contraction
from $\ell_2({\mathbb Z})$ into itself (the
system is then called {\sl dissipative}) if and
only if the function $\widehat{h}$ is analytic
and contractive in the open unit disk (such
functions are called {\sl Schur functions}), or
equivalently, if and only if the operator of
multiplication by $\widehat{h}$ is a contraction
from the Hardy space of the open unit disk
${\mathbf H}_2({\mathbb D})$ into itself. This
allows to resort to all the
tools of Schur analysis to study such systems; see for instance \cite{dd1}, 
\cite{MR2222523}, \cite{goh1}.\\

Note that $\widehat{h}$ is called in system theory the
{\sl transfer function} of the system. In certain fields
(e.g. engineering)
it is defined with $\zeta^{-1}$ instead of $\zeta$.\\

The system \eqref{patricia_kaas} commutes with
the shift operators
 $S$
\begin{equation}
\label{solferino} S(x_j)=(x_{j+1})
\end{equation}
defined in the input and output spaces. In
system theory terminology, it is called {\sl
time-invariant}. In fact, every time-invariant
linear bounded system from $\ell_2({\mathbb Z})$
into $\ell_2({\mathbb Z})$ is of this form. The
proof of this well known result is recalled in
the sequel; see
 STEP 1 in the proof of Theorem \ref{KSU_080909}.
Such a characterization does not hold when one
considers $\ell_\infty({\mathbb Z})$ instead of 
$\ell_2({\mathbb Z})$, as we will explain below.\\

The notion of Schur function and the associated
system theory interpretations, have been extended
in a number of directions, well beyond
the time-invariant case. We now discuss some of
them. First, when considering the time-varying
case (that is, when $h_{n-m}$ in
\eqref{patricia_kaas} is replaced by $h_{n,m}$),
an approach originating with the work of
Deprettere and Dewilde, see \cite{dede},
\cite{dede2}, consists of replacing the complex
numbers by diagonal operators. In the later
works \cite{MR92g:94002}, \cite{MR93b:47027},
the Hardy space is replaced by the Hilbert space
of upper-triangular operators of Hilbert-Schmidt
class, and and Schur functions by upper
triangular contractions. This allows, with an
appropriate definition of point evaluation of an
operator on a diagonal, to extend much of the
function theory of the open unit disk, to the
case of upper triangular operators, and hence to apply the
results to time-varying systems. See
\cite{DD-ot56}, \cite{bgk-ot56-1},
\cite{MR1911847}, \cite{MR1704663} for a sample
of papers, and \cite{MR99g:93001} for applications of this calculus on
diagonals.\\

Among other directions of research and extensions
we mention the case of multi-indexed systems and
their connections to several complex variables,
see for instance \cite{MR2184571}, and the
non-commutative case, see for instance
\cite{MR2129642}.\\

In all the directions outlined above, there is
no randomness in the system itself, although the
input (and hence the output) may be a sequence
of random variables. In the present work it is a
different kind of extension which we consider,
allowing the $h_n$ in \eqref{patricia_kaas}
to be random variables. We use white noise space
analysis, which has been introduced in 1975 by T.
Hida, see \cite{MR0451429}, and the monographs
\cite{MR1244577}, \cite{MR1408433} and
\cite{MR1387829}. White noise analysis allows to
translate problems from the stochastic context
into problems involving analytic functions in a
countable number of variables in the Fock space,
or in spaces of distributions which contain the
Fock space, in particular in the Kondratiev
space. The Wick product is a generalization for
random variables in the Kondratiev space of the
pointwise product, and reduces to the pointwise
product when at least one of the factors is nonrandom. It
became very useful when stochastic calculus with
respect to the fractional Brownian motion, and
more generally with respect to processes which
are not necessarily semi--martingales, began to
be considered; see \cite{MR1408433},
\cite{MR1801485}, \cite{MR1741154},
\cite{duncan1}.\\

Obviously, a Gaussian input into a linear system
with nonrandom coefficients, will result in a
Gaussian output. Here, we aim to model linear
Gaussian input-output relations when the
underlying linear system is random. While indeed
a Gaussian input into a linear system with random
coefficients cannot be expected to result in a
Gaussian output, we will use the white noise
space setting and replace the pointwise product
by the Wick product,
 enabling Gaussian
input-output relations when the underlying
system has random coefficients. This has the
advantage of preserving the Gaussian
input-output relation, while allowing
uncertainty in the form of randomness in the
linear system under study. Thus, we replace
\eqref{patricia_kaas} by
\begin{equation}
\label{system}
y_n= \sum_{m\in{\mathbb Z}}
h_{n-m}\lozenge u_m,\quad n\in\mathbb Z,
\end{equation}
where the $y_n$, $u_n$ and $h_n$ are now random
variables in the Kondratiev space (or more
precisely, in some Hilbert subspace of it), and
where $\lozenge$ denote the Wick product. We also
consider the causal case, where now
\begin{equation}
\label{system1}
y_n= \sum_{m=0}^n
h_{n-m}\lozenge u_m,\quad n=0,1,2,\ldots
\end{equation}

The proposed setting can be used to model uncertainty of an otherwise 
deterministic linear
time-invariant system, a system that maintains Gaussian input-output
relation, by a {\sl random} uncertainty in the impulse
response. This is known as the Bayesian embedding approach, by which
the study of a nonrandom uncertainty is carried out through an
associated probabilistic analysis; see e.g. \cite{MR690220}, \cite{MR1931659}.
We now turn to the content of the paper, and first recall three
stability theorems, namely, Theorems
\ref{celine_dion}, \ref{l1l2} and
\ref{avenue_parmentier}. The main aim of the
paper is to develop counterpart of these (and of
some other) theorems in the stochastic setting,
as is explained below.\\

Consider a linear discrete
time system of the form \eqref{patricia_kaas}.
Various notions of stability can be
assigned to such a system; in this work we will
focus in particular on BIBO stability (bounded
input bounded output), $\ell_1$-$\ell_2$
stability and the case of dissipative systems.
With BIBO stability in mind, the following result
is well known;  see for instance \cite[p.
177]{MR569473}.
\begin{Tm}
\label{celine_dion} There is a $M>0$ such that
the sums \( \sum_{m\in{\mathbb Z}} h_{n-m}u_m\)
converge absolutely for all $(u_m)\in\ell_\infty
({\mathbb Z})$, all $n\in{\mathbb Z}$, and
\begin{equation}
\label{bibo1} \sup_{n\in{\mathbb Z}}|y_n|\le
M\sup_{n\in{\mathbb Z}} |u_n|
\end{equation}
if and only if
\begin{equation}
\label{bibo} \sum_{n\in{\mathbb Z}} |h_n|\le M.
\end{equation}
\end{Tm}

Condition \eqref{bibo} means that the $\mathcal
Z$--transform \eqref{Ztr}
of the impulse response is in the Wiener algebra
${\mathscr W}$; $\widehat{h}$ is in particular
continuous on the unit circle, but it need not be
defined in general for $|\zeta|\not =1$. In the
case of a causal system, and when $h$ is
rational, \eqref{bibo} can be given a nicer
interpretation.
%
%
Recall first that a rational function which has
no pole on the unit circle is in the Wiener
algebra; its Taylor coefficients at the origin
(if the function is assumed analytic in a
neighborhood of the origin) will not, in
general, be equal to the Fourier coefficients of
its expansion as an element in $\mathscr W$.
They will be the same when all the $h_n=0$ for
$n$ negative, that is, when the system is
causal. Still for causal systems, condition
\eqref{bibo} means that the function
$\widehat{h}$ is analytic in the open unit disk
${\mathbb D}$, and continuous on the closed unit
disk. Thus, when $\widehat{h}$ is assumed rational, it
belongs to the subalgebra ${\mathscr W}_+$ of
$\mathscr W$. We also note that in this case,
$\widehat{h}$ has no pole in the closed unit disk. We
refer to \cite{ggk2} for these facts and for
more information on the
Wiener algebra.\\

The system \eqref{patricia_kaas} defines a linear
bounded operator from $\ell_\infty({\mathbb Z})$
into $\ell_\infty({\mathbb Z})$ which commutes
with the bilateral $S$ defined by
\eqref{solferino}; this last property expresses
the {\sl time-invariance} of the system. We
consider linear time-invariant systems whose
input-output relation is given in the form of a
convolution. We note however that not all linear
time-invariant systems from
$\ell_\infty({\mathbb Z})$ into itself are given
by a convolution. For instance,
define $X_0\in(\ell_\infty({\mathbb Z}))^*$ by
\[
X_0(x)=\limsup_{n\rightarrow\infty}\frac{
\sum_{j=-n}^{j=n} x_j}{2n+1}.
\]
The operator
\begin{equation}
\label{beilen} X(x)=(X_0S^{j*}(x))_{j\in{\mathbb Z}},
\end{equation}
defines a bounded linear operator from
$\ell_\infty({\mathbb Z})$ into itself. which
commute with $S$, but cannot be expressed by a convolution.\\

We denote
\[
{\mathbb N}=\left\{1,2,\ldots\right\}\quad {\rm
and}\quad {\mathbb N}_0={\mathbb N}\cup\left\{0
\right\}.
\]
To discuss $\ell_1$-$\ell_2$ stable and
dissipative systems, it is easier to consider
signals indexed by ${\mathbb N}_0$. We recall
that the Hardy space ${\mathbf H}_2({\mathbb
D})$ is the space of power series
\[
f(\zeta)=\sum_{n=0}^\infty f_n\zeta^n,\quad
f_n\in{\mathbb C},\]
with norm
\[
\|f\|_{{\mathbf H}_2({\mathbb
D})}=(\sum_{n=0}^\infty|f_n|^2)^{1/2}.
\]
The following definition and result are well
known:

\begin{Dn}
The system \eqref{patricia_kaas} will be called
$\ell_1$-$\ell_2$ bounded if there exists a
$M<\infty$ such that
\[
\left(\sum_{n=0}^\infty
|y_n|^2\right)^{1/2}\le
M\sum_{n=0}^\infty|u_n|.
\]
\end{Dn}

For the matrix-valued version of the theorem below, see
\cite[Theorem 5.1]{Aleb2}.
\begin{Tm}
\label{l1l2} The system \eqref{system} is
$\ell_1$-$\ell_2$ bounded if and only if its
transfer function is in the Hardy space
${\mathbf H}_2({\mathbb D})$.
\end{Tm}

The $\ell_2({\mathbb N}_0)$ norms of the input
and output sequences are a measure of the energy
of the signals, and play an important role in
system theory; see \cite{MR2002b:47144} for a survey.
The system \eqref{patricia_kaas} is called
dissipative if the $\ell_2({\mathbb N}_0)$ norm
of the output is always less or equal to the
$\ell_2({\mathbb N}_0)$-norm of the input. The
following result characterizes systems of the
form \eqref{patricia_kaas} which are dissipative.

\begin{Tm}
\label{avenue_parmentier} A linear system is
time-invariant, causal and dissipative if and
only if it is of the form \eqref{patricia_kaas}
with a transfer function which is analytic and
contractive in the open unit disk.
\end{Tm}

In other words, the system has a transfer
function which is a Schur function.
Equivalently, the lower triangular Toeplitz
operator
\begin{equation}
\label{toep_080909}
\begin{pmatrix}h_0&  0    &\cdot&\cdots\\
               h_1&h_0&\cdot&\cdots\\
               h_2&h_1&h_0&\cdots\\
 \cdot&\cdot&\cdot&\cdots
\end{pmatrix}
\end{equation}
is a contraction from
$\ell_2({\mathbb N}_0)$ into itself.\\

We note that, through the $\mathcal Z$-transform,
the space $\ell_2({\mathbb N}_0)$ is unitarily
mapped onto
${\mathbf H}_2({\mathbb D})$, that is, onto the
reproducing kernel Hilbert space with
reproducing kernel
\begin{equation}
\label{hardy}
K(\zeta,\nu)=\frac{1}{1-\zeta\nu^*},\quad
\zeta,\nu\in{\mathbb D}.
\end{equation}
Therefore, $\widehat{h}$ is the transfer function
of a dissipative system if and only if the
operator of multiplication by $\widehat{h}$ is a
contraction from ${\mathbf H}_2({\mathbb D})$
into itself, or, equivalently, if and only if
the kernel
\begin{equation}
\frac{1-\widehat{h}(\zeta)
\widehat{h}(\nu)^*}{1-\zeta\nu^*}
\label{Manhattan_090909}
\end{equation}
is positive in the open unit disk. The
associated reproducing kernel Hilbert spaces
were introduced and studied by de Branges and
Rovnyak, also in the operator-valued case; see
\cite{dbr1}, \cite{dbr2}. We will use some of
their results in the sequel; see Theorem
\ref{ma}.
\\

We wish to extend the notion of transfer function
so as to include a random aspect, and present
counterparts of the three theorems mentioned
above in a random systems setting. We will also
consider the continuous case for BIBO stability,
dissipative systems, and ${\mathbf
L}_2$-${\mathbf L}_\infty$ bounded systems.
Using the white noise setting and the Hermite
transform this random aspect will be expressed
by dependence on a countable number of
independent complex variables $z_1,z_2,\ldots$.
This is reviewed in the next section. First we
need to recall some of the relevant notations at
this stage. In the following expressions, $\ell$
denotes the set of finite sequences of integers
indexed by ${\mathbb N}$, that is, the set of
sequences $(n_1,n_2,\ldots)$ where $n_j=0$ for
all large enough $j$, and we set (see
\cite[(2.3.8) p. 29]{MR1408433})
\[
(2{\mathbb N})^\alpha=\prod_{j\in{\mathbb N}}
(2j)^{\alpha_j}.\] The product is meaningful
since $\alpha_j=0$ for all but for a finite
number of indices $j$.\\

In our approach, we replace the kernel
\eqref{hardy} by a kernel of the form
\[
\frac{K_k(z,w)}{1-\zeta\nu^*}
\]
where $z=(z_1,z_2,\ldots )$ and
$w=(w_1,w_2,\ldots)$ belong to the infinite
dimensional neighborhood
\[
{\mathbb K}_{k}=
\left\{z=(z_1,z_2,\ldots)\in{\mathbb
  C}^{\mathbb N}\,
:\, \sum_{\alpha\in\ell} |z|^{2\alpha}(2{\mathbb
N})^{k\alpha}<\infty \right\}
\]
of the origin in ${\mathbb C}^{\mathbb N}$ (see
\cite[Definition 2.6.4 p. 59]{MR1408433}), and
where
\begin{equation}
\label{sdf}
K_k(z,w)=\sum_{\alpha\in\ell}z^\alpha(w^*)^\alpha
(2{\mathbb N})^{k\alpha},
\end{equation}
with the use of the multi-index notation:
\[
z^\alpha=z_1^{\alpha_1}\cdot
z_2^{\alpha_2}\cdots,\quad{\rm where}\quad
\alpha=(\alpha_1,\alpha_2,\ldots).
\]
We denote by ${\mathcal H}(K_k)$ the reproducing
kernel Hilbert space with reproducing kernel
$K_k(z,w)$. Elements of ${\mathcal H}(K_k)$ are
analytic in $z=(z_1,z_2,\ldots)$ in the set
${\mathbb K}_k$. To take into account
randomness, we replace the Hardy space ${\mathbf
H}_2$ by ${\mathbf H}_2 \otimes {\mathcal
H}(K_k)$. Note that
\begin{equation}
\label{la_muette} {\mathbf H}_2 \otimes {\mathcal
H}(K_k) =
\left\{ f(\zeta)=\sum_{n=0}^\infty \zeta^n
f_n\,{\rm with}\quad
 f_n\in{\mathcal H}(K_k)\right\},
\end{equation}
with norm
\[
\|f\|^2_{{\mathbf H}_2 \otimes {\mathcal
H}(K_k)} =\sum_{n=0}^\infty\|f_n\|^2_{{\mathcal
H}(K_k)}.
\]
A transfer function will thus be a function
${\mathscr H}(\zeta,z)$ of the form
\begin{equation}
{\mathscr H}(\zeta,z)=\sum_{n=0}^\infty \zeta^n
h_n(z), \label{Mouton_Duvernet_ligne_4}
\end{equation}
where now the $h_n\in{\mathcal H}(K_k)$. In
\eqref{Mouton_Duvernet_ligne_4}, the nonrandom
part of the transfer function is
\begin{equation}
\label{Porte} {\mathscr H}(\zeta,
0)=\sum_{n=0}^\infty \zeta^n h_n(0).
\end{equation}

As we explain in the sequel, another possible
generalization is to write
\begin{equation}
\label{madeleine} {\mathscr
H}(\zeta,z)=\widehat{h}(\zeta) +\sum_{\alpha\in
\ell} z^\alpha c_\alpha(\zeta).
\end{equation}

The nonrandom part in \eqref{madeleine} is
\begin{equation}
\label{paris} {\mathscr H}(\zeta,
0)=\widehat{h}(\zeta),
\end{equation}
corresponding to $\alpha=0$.\\


The transfer function ${\mathscr H}(\zeta,z)$ of
a dissipative random system will be
characterized by the positivity of the kernel
\begin{equation}
(1-{\mathscr H}(\zeta,z) {\mathscr
H}(\nu,w)^*)\frac{K_k(z,w)} {(1-\zeta\nu^*)}
\label{chaud_les_marrons}
\end{equation}
in ${\mathbb D}\times{\mathbb K}_k$, that is, ${
\mathscr H}(\zeta, z)$ is a contractive
multiplicator of the reproducing kernel Hilbert
space ${\mathbf H}_2({\mathbb D})\otimes
{\mathcal H}(K_k)$ with reproducing kernel
\begin{equation}
\label{republique111}
\frac{K_k(\zeta,z)}{1-\zeta\nu^*}.
\end{equation}
A similar interpretation will hold in the
continuous case, where the open unit disk is now
replaced by the open upper
half-plane.\\

The outline of the paper is as follows. The
paper consists of eight sections besides the
introduction. In the second section, we review
the white space noise setting, and define
analogs of the linear systems
\eqref{patricia_kaas} and of their continuous
time versions when we allow randomness both in
the impulse response $(h_m)$ and in the inputs
$(u_m)$. In Section 3 we prove the counterparts
of Theorems \ref{celine_dion}. The cases of
Theorems \ref{l1l2} and \ref{avenue_parmentier}
are considered in Section 4 and 5 respectively.
The next three sections are devoted to the
continuous time case. In Section 6 we prove the
operator-valued version of the
Bochner-Chandrashekaran theorem. In Section 7 we
consider the analog of BIBO continuous systems
and in Section 8 the case of continuous
dissipative systems. In the last section we
consider the case of ${\mathbf L}_2$-${\mathbf
L}_\infty$
stability.\\

The theory of random linear systems is yet not
well established, and it appears that the
setting presented in this paper is the first
which permits to handle, in an efficient way, the
case where randomness is allowed in the impulse
response.\\

We denote by
\[
Uf(t)=\frac{1}{\sqrt{2\pi}}\int_{\mathbb R}
f(x)e^{ixt}dx
\]
the Fourier transform, where $f\in{\mathbf
L}_2({\mathbb R})$, or, more generally, belongs
to ${\mathbf L}_2({\mathbb R})\otimes{\mathcal
H}$ for some Hilbert space ${\mathcal H}$.
For a Hilbert space ${\mathcal H}$, the symbol
${\mathbf L}({\mathcal H})$ denotes the set of
bounded operators from ${\mathcal H}$ into
itself.

\section{The white noise space setting}
\setcounter{equation}{0}

We now briefly recall the definitions of the
white noise space and of the Kondratiev space.
We refer to \cite{MR1244577},
\cite{MR1408433},\cite{MR1387829}, \cite{bosw}
and
\cite{eh} for more information.\\

The function
\[
K(s_1-s_2)= \exp(-\|s_1-s_2\|_{{\bf
L}_2({\mathbb R})}^2/2)
\]
is positive on the Schwartz space ${\mathcal S}$
of rapidly decreasing, infinitely differentiable
functions. Since ${\mathcal S}$ is a nuclear
space, the Bochner-Minlos theorem (see \cite[Th\'eor\`me 2 p. 342]{MR35:7123})
ensures the
existence of a probability measure $P$ on the
Borel sets ${\mathcal
  F}$ of
$\Omega={\mathcal S}^\prime$ such that
\[
\int_{\Omega}e^{i\langle s,\w \rangle_{{\bf
L}_2({\mathbb R})}}dP(\w) =\exp(-\|s \|_{{\bf
L}_2({\mathbb R})}^2/2).
\]
The white noise space is ${\mathcal W}={\mathbf
L}_2(\Omega, {\mathcal F}, P)$.
 An orthogonal basis of the white
noise space is given by the Hermite functions
$(H_\alpha)_{\alpha\in\ell}$. These functions
are computed in terms of the Hermite polynomials,
and we refer to  \cite[Definition 2.2.1 p.
19]{MR1408433} for their definition. Every
element in ${\mathcal W}$ can be written as
\begin{equation}
\label{fps} F(\w)=\sum_{\alpha\in\ell}c_\alpha
H_\alpha(\w),\quad c_\alpha\in{\mathbb C},
\quad{\rm with}\quad \|F\|^2_{{\mathcal
W}}=\sum_{\alpha\in\ell} |c_\alpha|^2\alpha!<\infty
\end{equation}
In general, one takes real $c_\alpha$. Here we
take complex coefficients, that is, we consider
the complexified real white noise space. A
similar remark holds for the spaces ${\mathcal
H}(K_k)$
defined below.\\

Let $z=(z_1,z_2,\ldots)\in{\mathbb C}^{\mathbb
N}$. The linear map which to $H_\alpha$
associates the polynomial ${\mathbf
I}(H_\alpha)=z^\alpha$ extends to a unitary map
between the white noise space and the real
reproducing kernel Hilbert space
with reproducing kernel $\exp{\langle z,
w\rangle_{\ell_2}}$. This space is called the
Fock space. The map ${\mathbf I}$ is called the
{\sl Hermite transform}. The Wick product is
defined through the Hermite functions by
\[
H_\alpha\lozenge H_\beta= H_{\alpha+\beta},\quad
\alpha,\beta\in\ell.
\]
It is in fact independent of the chosen basis in
the white noise space; see \cite[Appendix D, pp.
209-215]{MR1408433}. In general, the Wick
product of two elements in the white noise space
need not be in the white noise space. The most
convenient space which is stable with respect to
the Wick product is the Kondratiev space
$S_{-1}$. To define $S_{-1}$ we first introduce
for $k\in{\mathbb N}=\left\{1,2,\ldots\right\}$
the Hilbert space ${\mathcal H}_{k}$ which
consist of series of the form \eqref{fps} such
that
\begin{equation}
\label{michelle}
\|f\|_{k}\stackrel{\rm def.}{=}
 \left(\sum_{\alpha\in\ell}|c_\alpha|^2
(2{\mathbb N})^{-k\alpha}\right)^{1/2}<\infty.
\end{equation}

Note that the the Hermite transform is a unitary mapping
from ${\mathcal H}_{k}$ onto the reproducing kernel
Hilbert space with reproducing kernel $K_k(z,w)$,
where $K_k$ is defined in \eqref{sdf}. The
Kondratiev space $S_{-1}$ is the inductive limit
of the spaces ${\mathcal H}_{k}$. We note that
when either one of the factors $f$ or $g$ in
$S_{-1}$ is nonrandom, the Wick product
$f\lozenge g$ reduces to the pointwise product
$fg$.\\

We will consider stochastic processes in the
series form
\begin{equation}
\label{republique} X(\tau,\w)=
\sum_{\alpha\in\ell}c_\alpha(\tau)H_\alpha(\w),
\end{equation}
where the $c_\alpha(\tau)$ are nonrandom
functions depending on a parameter $\tau$. We
require that the series \eqref{republique}
belongs to the Kondratiev space $S_{-1}$ for
every value of $\tau$. Here $\tau$ belongs to
the integers or to the real numbers. We note that
other choices of indices are possible (for
example, ${\mathbb Z}^2$, or the case of the
vertices of a binary tree). We also note that a
case of special interest arises when one applies the Hermite
transform and the Laplace transform (for the
continuous time case) or the $\mathcal
Z$-transform (for the discrete time case), to
obtain a function which depends on a finite
number of variables $z_i$, which, in addition, is a rational
function in these variables, that is, is a function of the
form

\begin{equation}
\label{toto} D(\zeta)+C(\zeta)(I_N-\sum_{k=1}^M
z_kA_k(\zeta))^{-1}(\sum_{k=1}^Mz_kB_k(\zeta)).
\end{equation}
Here, $\zeta$ denotes the variable corresponding either to the
Laplace transform, or the $\mathcal
Z$-transform. The functions of $\zeta$ in that
expression are also assumed rational. See
\cite{alpay-2008}. Other type of realizations
are possible; see Theorem \ref{ma} below.

\section{BIBO stable linear
discrete time stochastic systems}
\setcounter{equation}{0}
Fix some integer $l>0$,  and let $k>l+1$. Consider
$h\in {\mathcal H}_{l}$ and $u\in {\mathcal
H}_{k}$. Then, V\r{a}ge's inequality (see
\cite[Proposition 3.3.2 p. 118]{MR1408433}) is
in the form
\begin{equation}
\label{vage} \|h\lozenge u\|_{k}\le
A(k-l)\|h\|_{l}\|u\|_{k},
\end{equation}
where
\begin{equation}
\label{vage111}
A(k-l)=\sum_{\alpha\in\ell}(2{\mathbb
N})^{(l-k)\alpha}.
\end{equation}
It is not a trivial fact that $A(k-l)$ is
finite; see \cite{MR1275367} and
\cite[Proposition 2.3.3 p. 31]{MR1408433} for a proof.\\

Inequality \eqref{vage} expresses the fact that
the multiplication operator
\[
T_h\,:\,u\mapsto h\lozenge u
\]
is a bounded map from the Hilbert space
${\mathcal H}_{k}$ into itself, and that its
operator norm $\|T_h\|_{{\rm op},l, k}$ satisfies
the inequality
\begin{equation}
\label{belmondo} 
\|T_h\|_{{\rm op},l,k}\le
A(k-l)\|h\|_l.
\end{equation}
The norm of $T_h$ depends on $k$ and $l$ and
will be in general different from $\|h\|_{l}$.
It implies in particular that ${\mathcal H}_l$
endowed with the norm
\[
\|h\|\stackrel{\rm def.}{=}\|T_h\|_{{\rm op},l,k}
\]
is a normed algebra.\\

To simplify the notation, we set
\[
\|T_h\|_{{\rm op},l,k}=\|T_h\|.
\]

\begin{Dn}
A random discrete time signal will be a sequence
$(u_n)$ indexed by ${\mathbb Z}$, of elements in
the Kondratiev space, such that there exists a
$k\in{\mathbb N}$ (depending on the signal) such
that
\[
u_n\in{\mathcal H}_{k},\quad  \forall
n\in{\mathbb Z}.
\]
\end{Dn}

We note that $k$ is imposed to be independent of
$n$.
\begin{Tm}
Let $k>l+1$ and let $(h_n)$ be a sequence of
elements in ${\mathcal H}_{l}$ indexed by
$\mathbb Z$.
Then \\
$(a)$ The sums \eqref{system}
\begin{equation*}
y_n= \sum_{m\in{\mathbb Z}} h_{n-m}\lozenge
u_m,\quad n\in\mathbb Z,
\end{equation*}
converge absolutely in ${\mathcal H}_k$ for all
inputs $(u_m)_{m\in{\mathbb Z}}$ such that
\(\sup_{m\in\mathbb Z}\|u_m\|_{k}<\infty\),
and\\
$(b)$ There exists an $M>0$ such that, for all such
inputs $(u_n)_{n\in{\mathbb Z}}$, it holds that
\begin{equation}
\label{220708} \sup_{n\in{\mathbb
Z}}\|y_n\|_{k}\le M\sup_{ n\in{\mathbb Z}}
\|u_n\|_{k}
\end{equation}
if and only if for all $v\in {\mathcal H}_k$ with
$\|v\|_{k}=1$ it holds that
\begin{equation}
\label{200708} \sum_{n\in{\mathbb Z}}
\|T_{h_n}^*(v)\|_{k}\le M.
\end{equation}
\label{joan_baez}
\end{Tm}

We note that, in the nonrandom case, where the
$h_n$ are (nonrandom) complex numbers, the Wick
product reduces to a pointwise product, and we
have systems of the form \eqref{patricia_kaas}.
Furthermore, in this case,
\[
\|T_{h_n}^*v\|_{k}=|h_n|\cdot \|v\|_{k},
\]
and we retrieve the well known BIBO stability
condition \eqref{bibo}, with Theorem
\ref{joan_baez} reduced to
Theorem \ref{celine_dion}.\\

Furthermore, we notice that condition
\[
\sum_{n\in{\mathbb Z}}\|T_{h_n}\|\le M
\]
on the norms of the operators $T_{h_n}$ implies
condition \eqref{200708}. Expressions of the
form
\[
\sum_{n\in\mathbb Z}\zeta^n T_{h_n}
\]
with
\[
\sum_{n\in\mathbb Z}\|T_{h_n}\|<\infty
\]
form an algebra, which appears to be the
counterpart of the classical Wiener algebra in
the present setting.\\

{\bf Proof of Theorem \ref{joan_baez}:} First
note that, in view of the restriction $k>l+1$,
each of the terms in \eqref{system} belongs to
${\mathcal H}_k$.\\

Assume that \eqref{220708} is in force. Then, for
every $n\in{\mathbb Z}$, taking sequences
$(u_m)$ which have only a finite number of non
zero entries, one has for any preassigned $n \in
{\mathbb Z}$,

\[
\begin{split}
M\sup_{m\in{\mathbb Z}}\|u_m\|_{k}&\ge
\|\sum_{m\in{\mathbb Z}} T_{h_{n-m}}\lozenge u_m
\|_{k}\\&=\sup_{\substack{v\in {\mathcal H}_k\\
\|v\|_{k}=1}}\langle \sum_{m\in{\mathbb Z}}
T_{h_{n-m}}\lozenge u_m,
v\rangle_{k}\\
&= \sup_{\substack{v\in {\mathcal H}_k\\
\|v\|_{k}=1}}\langle
\sum_{\substack{m\in{\mathbb Z}\\
  u_m\not =0}}
u_m, T_{h_{n-m}}^*(v)\rangle_{k}.
\end{split}
\]
The special choice
\[
u_m=\begin{cases}
\frac{T_{h_{n-m}}^*(v)}{\|T_{h_{n-m}}^*(v)
\|_{k}},\quad{\rm if}\quad
{\|T_{h_{n-m}}^*(v) \|_{k}}\not =0\\
\hspace{7mm}0,\hspace{1.cm}\quad{\rm if}\quad
{\|T_{h_{n-m}}^*(v) \|_{k}} =0
\end{cases}
\]
leads to
\[
\sum_{\substack{m\in{\mathbb Z}\\
  u_m\not =0}}
\|T_{h_{n-m}}^*(v)\|_{k}\le M.
\]
This expression stays the same also for the indices $m$ such that
$u_m\not=0$  since $u_m\not=0$ if and only if 
$\|T_{h_{n-m}}^*(v) \|_{k}\not=0$. Finally, since the right handside
of the above inequality is independent of the support of $(u_m)$ we
get the result.\\

Conversely, assume that \eqref{200708} is in
force. Then, still for sequences $(u_m)$ with
only a finite number of non zero entries, we have
\[
\begin{split}
\|\sum_{\substack{m\in{\mathbb Z}\\
  u_m\not =0}}
T_{h_{n-m}}u_m\|_{k}&= \sup_{\substack{v\in
{\mathcal H}_k\\ \|v\|_{k}=1}}
\langle \sum_{\substack{m\in{\mathbb Z}\\
  u_m\not =0}}
u_m,
T_{h_{n-m}}^*(v)\rangle_{k}\\
&\le\sup_{\substack{v\in {\mathcal H}_k\\
\|v\|_{k}=1}}
\sum_{\substack{m\in{\mathbb Z}\\
  u_m\not =0}}
\|u_m\|_{k}\|T_{h_{n-m}}^*v\|_{k}\\
&\le \sup_{\substack{m\in{\mathbb Z}\\
  u_m\not =0}}
\|u_m\|_{k}\cdot
\sup_{\substack{v\in {\mathcal H}_k\\
\|v\|_{k}=1}} \left( \sum_{m\in{\mathbb Z}}
\|T_{h_m}^*(v)\|_{k}\right)\\
&\le M \sup_{\substack{m\in{\mathbb Z}\\
  u_m\not =0}}\|u_m\|_{k}.
\end{split}
\]
Assume now that the sum
\[
\sum_{m\in{\mathbb Z}}
T_{h_{n-m}}u_m
\]
converges absolutely in ${\mathcal H}_{k}$.
The result is then obtained by continuity by 
considering partial finite sums.
\mbox{}\qed\mbox{}\\

Applying the Hermite transform to equations
\eqref{system}, that is to
\[
y_n=\sum_{m\in{\mathbb Z}}h_{n-m}\lozenge
u_m,\quad n\in{\mathbb Z},
\]
leads to the following: let
\[
y_n(\w)=\sum_{\alpha\in\ell}y_\alpha(n)
H_\alpha(\w)\quad {\rm and}\quad h_n(\w)=
\sum_{\alpha\in\ell}h_\alpha(n)H_\alpha(\w),\]
where the coefficients $y_\alpha(n)$ and
$h_\alpha(n)$ are nonrandom complex numbers.
Then,
\begin{equation*}
y_n=\sum_{\alpha\in\ell} H_\alpha(\w)
\left(\sum_{m\in{\mathbb Z}}\sum_{\beta\le
\alpha} h_{\alpha-\beta}(n-m)u_\beta(m)\right),
\quad n\in{\mathbb Z},
\end{equation*}
that is, after taking the Hermite transform
\begin{equation*}
{\mathbf I}(y_n)=\sum_{\alpha\in\ell} z^\alpha
\left(\sum_{m\in{\mathbb Z}}\sum_{\beta\le
\alpha} h_{\alpha-\beta}(n-m)u_\beta(m)\right),
\quad n\in{\mathbb Z},
\end{equation*}
and hence,
\begin{equation}
y_\alpha(n)=\sum_{m\in{\mathbb Z}}\sum_{\beta\le
\alpha} h_{\alpha-\beta}(n-m)u_\beta(m),\quad
n\in{\mathbb Z}. \label{paris_texas}
\end{equation}
Equation \eqref{paris_texas} exhibits two
convolutions: the first is with
respect to the index in $\ell$, which is related
to the stochastic aspect of the system; the
second is with respect to the time
variable.
\\

The ${\mathcal Z}$  transform (denoted by
$\widehat{y}$, with variable $\zeta$) then
leads to

\[
\begin{split}
\widehat{y}(\zeta,z)&\stackrel{\rm
def.}{=}\sum_{n\in{\mathbb Z}}{\mathbf
I}(y_n)\zeta^n\\
& =\sum_{\alpha\in\ell} z^\alpha
\sum_{\beta\le\alpha}\widehat{h}_{\alpha-\beta}
(\zeta)\widehat{u}_\beta(\zeta)\\
&=\left(\sum_{\alpha\in\ell}
z^\alpha\widehat{h}_\alpha(\zeta)\right)
\left(\sum_{\alpha\in\ell}
z^\alpha\widehat{u}_\alpha(\zeta)\right).
\end{split}
\]
\begin{Dn}
The function
\[
{\mathscr H}(\zeta, z) =\sum_{\alpha\in\ell}
z^\alpha\widehat{h}_\alpha(\zeta)
=\sum_{n\in{\mathbb Z}}\zeta^n({\mathbf
I}(h_n))(z)
\]
is the called generalized transfer function of the
system.
\end{Dn}

When all $\widehat{h}_\alpha(\zeta)=0$
for $\alpha\not =0$, we retrieve the classical
notion of the transfer function. We can thus
define a hierarchy of systems, depending on the
properties of the function ${\mathscr H}(\zeta,
z)$. The rational case will be when the function
${\mathscr H}(\zeta, z)$ is of the form
\eqref{toto}. Another case of interest would be
the isospectral case, when the function
$A(\zeta)$ in \eqref{toto} does not depend on the variable $\zeta$.\\

\section{$\ell_1$-$\ell_2$
stable random systems}
\setcounter{equation}{0}
The analog of Theorem \ref{l1l2} is the
following:
\begin{Tm}
Let $l>k+1$ and assume that in the system \eqref{system1} 
$h_n\in {\mathcal H}_l$.
Then there exists $M>0$ such that
\[
\left(\sum_{n=0}^\infty \|y_n\|_k^2\right)^{1/2}
\le M\sum_{n=0}^\infty \|u_n\|_k
\]
for all inputs $(u_n)$ such that the right handside of the above equation 
is finite,
if and only if its
transfer function belongs to ${\mathbf
H}_2({\mathbb D})\otimes {\mathcal H}(K_l)$, i.e.
if and only if there exists a number $c>0$ such that the
kernel
\begin{equation}
\label{trocadero}
\frac{K_l(z,w)}{1-\zeta\nu^*}-c{\mathscr
H}(\zeta,z){\mathscr H}(\nu,w)^*
\end{equation}
is positive in ${\mathbb D}\times {\mathbb K}_l$.
\label{salomon}
\end{Tm}

The system \eqref{system} is then called $\ell_1$-$\ell_2$ bounded. 
In the proof of the theorem, we use the already mentioned fact that the spaces
${\mathcal H}_k$ and ${\mathcal H}(K_k)$ are unitarily equivalent via
the Hermite transform. This allows us to make use of V\r{a}ge's
inequality in the spaces ${\mathcal H}(K_k)$.\\

{\bf Proof of Theorem \eqref{salomon}:}  We first remark that the
equivalence of \eqref{trocadero} with the
condition ${\mathscr H}\in {\mathbf H}_2(
{\mathbb D})\otimes {\mathcal H}(K_k)$ directly
follows the characterization of the elements of a
reproducing kernel space. We proceed in three
steps. Recall that $A(k-l)$
has been defined by \eqref{vage111}.\\

STEP 1: {\sl Assume that \eqref{trocadero} 
holds, and let $u\in{\mathcal H}(K_k)$. Then,
\begin{equation}
\label{ineq111} \|{\mathscr H}x\|_{{\mathbf H}_2(
{\mathbb D})\otimes {\mathcal H}(K_k)} \le
A(k-\ell) \|{\mathscr H} \|_{{\mathbf H}_2(
{\mathbb D})\otimes {\mathcal
H}(K_l)}\cdot\|u\|_{{\mathcal H}(K_k)}. 
\end{equation}
}

Indeed, let $ {\mathscr H}(\zeta,
z)=\sum_{n=0}^\infty \zeta^n h_n(z)$, with
$h_n\in {\mathcal H}(K_\ell)$. Then,
\[
\begin{split}
\|{\mathscr H}u\|^2_{{\mathbf H}_2( {\mathbb
D})\otimes {\mathcal H}(K_k)}&=\sum_{n=0}^\infty
\|h_nu\|^2_k\\
&\le A(k-l)^2\sum_{n=0}^\infty
\|h_n\|_{{\mathcal H}(K_l)}^2\|u\|^2_{{\mathcal H}(K_k)}\\
&=A(k-l)^2\|\mathscr H \|_{{\mathbf H}_2(
{\mathbb D})\otimes{\mathcal
H}(K_l)}\cdot\|u\|^2_{{\mathcal H}(K_k)}.
\end{split}
\]

STEP 2: {\sl Assume that \eqref{trocadero} holds.
Then, the system is $\ell_1$-$\ell_2$ bounded.}\\

We first note that, from the definition of the
norm in ${{\mathbf H}_2( {\mathbb D})\otimes
{\mathcal H}(K_l)}$, the operator $M_\zeta$ of
multiplication by $\zeta$ is an isometry from
${{\mathbf H}_2( {\mathbb D})\otimes {\mathcal
H}(K_l)}$ into itself:
\begin{equation}
\label{isometrie} \|M_\zeta F\|_{{\mathbf H}_2(
{\mathbb D})\otimes {\mathcal H}(K_l)}=
\|F\|_{{\mathbf H}_2( {\mathbb D})\otimes
{\mathcal H}(K_l)},\quad \forall F\in{{\mathbf
H}_2( {\mathbb D})\otimes {\mathcal H}(K_l)}.
\end{equation}

Let $\widehat{u}(\zeta, z)=\sum_{n=0}^\infty \zeta^n
u_n(z)$, where the $u_n\in{\mathcal H}(K_k)$. Then,
\[
\widehat{y}(\zeta,z)=\sum_{n=0}^\infty \mathscr H(\zeta,
z)\zeta^nu_n(z),
\]
where the convergence is in ${{\mathbf H}_2( {\mathbb
D})\otimes {\mathcal H}(K_l)}$, as will be
justified by equation \eqref{conv} below. So,
using \eqref{ineq111} and \eqref{isometrie}, we
have, with $\widehat{y}(\zeta,z)=\sum_{n=0}^\infty
\zeta^ny_n(z)$,
\begin{equation}
\label{conv}
\begin{split}
\left(\sum_{n=0}^\infty\|y_n\|^2_{{\mathcal
H}(K_k)}\right)^{1/2}&=
\|\widehat{y}\|_{{\mathbf H}_2(
{\mathbb D})\otimes {\mathcal H}(K_k)}\\
&=\sum_{n=0}^\infty\|{\mathscr H}M_\zeta^n u_n\|_{{\mathbf
H}_2( {\mathbb D})\otimes{\mathcal H}(K_l)}\\
&\le\sum_{n=0}^\infty\|{\mathscr H}u_n\|_{{\mathbf
H}_2( {\mathbb D})\otimes{\mathcal H}(K_l)}\\
&\le
A(k-l) \sum_{n=0}^\infty\|\mathscr H\|_{{\mathbf
H}_2( {\mathbb
D})\otimes {\mathcal H}(K_l)}\|u_n\|_{{\mathcal H}(K_k)}\\
&=M\left(\sum_{n=0}^\infty\|u_n\|_{{\mathcal H}(K_k)}\right),
\end{split}
\end{equation}
with $M=A(k-l)\|\mathscr H\|_{{\mathbf H}_2(
{\mathbb
D})\otimes {\mathcal H}(K_l)}$.\\

STEP 3: {\sl We complete the proof.}\\

Assume the system $\ell_1$-$\ell_2$ bounded.
Since $1\in{\mathcal H}(K_k)$ we have that \(
{\mathscr H}(\zeta,z) \in{{\mathbf H}_2( {\mathbb
D})\otimes {\mathcal H}(K_k)} \), and we can use
the previous step to conclude the proof.
\mbox{}\qed\mbox{}\\

When the system is not random, that is, when
$\mathscr H$ is only a function of $\zeta$, the
positivity of the kernel \eqref{trocadero} is
equivalent to the fact that ${\mathscr H}\in
{\mathbf H}_2({\mathbb D})$. Indeed, the tensor
product space ${\mathbf H}_2({\mathbb D})\otimes
{\mathcal H}(K_k)$ can be seen not only as the
representation \eqref{la_muette}, but also as the
space of functions of the form
\[
F(\zeta,z)=\sum_{\alpha\in\ell}z^\alpha f_\alpha,
\quad f_\alpha\in{\mathbf H}_2({\mathbb D}),\]
with the norm
\[
\|F\|_{{\mathbf H}_2({\mathbb D})\otimes
{\mathcal H}(K_k)}^2=\sum_{\alpha\in\ell}\frac{
\|f_\alpha\|^2_{{\mathbf H}_2({\mathbb
D})}}{(2{\mathbb N})^{k\alpha}}.
\]

\section{Dissipative discrete time random systems}
\setcounter{equation}{0}
We denote by $M_{\zeta}$ the operator of
multiplication by $\zeta$ and by $M_{z_j}$ the
operator of multiplication by $z_j$ (in both
cases, the domain and range of these operators
is given below). In the next theorem, the first
condition expresses the time-invariance of the
system. The second condition is needed to ensure
that we get a multiplication operator. It can be
interpreted as the property of invariance with respect to
randomness. In the following statement, recall
that the kernel $K_k$ has been defined by
\eqref{sdf}.
\begin{Tm}
Let $k\in{\mathbb N}$. The operators $M_{z_j}$
are bounded from ${\mathcal H}(K_k)$ into
itself. A linear operator from ${\mathbf
H}_2({\mathbb D})\otimes {\mathcal H}(K_k)$ into
itself is contractive and such that
\begin{equation}
\label{luxembourg}
\begin{split}
T(M_\zeta f)&=M_\zeta Tf\\
T(M_{z_j}f)&=M_{z_j}Tf
\end{split}
\end{equation}
if and only if it is of the form
\[
(Tf)(\zeta,z)={\mathscr S}(\zeta, z)f(\zeta,z)
\]
where $\mathscr S$ is such that the kernel
\eqref{chaud_les_marrons}:
\[ (1-{\mathscr
S}(\zeta,z) {\mathscr S}(\nu,w)^*)\frac{K_k(z,w)}
{(1-\zeta\nu^*)}
\]
is positive in ${\mathbb D}\times {\mathbb K}_k$.
\label{KSU_080909}
\end{Tm}
{\bf Proof:} We divide the proofs into several
steps.\\

STEP 1: {\sl Let ${\mathcal H}$ be a Hilbert
space, and let $T$ be a bounded operator from
${\mathbf H}_2({\mathbb D})\otimes{\mathcal H}$
into itself which commutes with multiplication by
$\zeta$. Then, $T$ is of the form
\[
Tf(\zeta)=S(\zeta)f(\zeta),\]
where $S$ is an ${\mathbf L}({\mathcal
H})$-valued function analytic and contractive in
the open unit
disk.}\\

This is a well known fact (see for instance
\cite[Lemma 1 p. 301]{MR0140931}), being the
discrete analog of the Bochner-Chandrasekharan
theorem (see \cite[Theorem 72, p.
144]{boch_chan} and Section \ref{rotterdam} below
for the latter). We briefly
review its proof for completeness. Let
$c\in{\mathcal H}$. We have
\[
(Tc)(\zeta)=\sum_{n=0}^\infty \zeta^n T_n(c),
\]
where the operators $T_n$ are readily seen to be
linear bounded operators from ${\mathcal H}$ into
itself. Furthermore,
\[
\|T_n\|\le \|T\|,\quad\forall n\in{\mathbb N}_0,
\]
since
\[
\|Tc\|^2_{{\mathbf H}_2({\mathbb
D})\otimes{\mathcal H}}=
\sum_{m=0}^\infty\|T_mc\|^2_{\mathcal H}\ge
\|T_nc\|^2_{{\mathcal H}},\quad\forall
c\in{\mathcal H}\quad{\rm and}\quad\forall n\in
{\mathbb N}_0.
\]
The ${\mathbf L}({\mathcal H})$--valued function
\begin{equation}
\label{ST}
S(\zeta)=\sum_{n=0}^\infty \zeta^n T_n
\end{equation}
is analytic in the open unit disk (we refer to
\cite[pp. 189--190]{MR58:12429a} for a review of
operator-valued analytic functions, and the
equivalence between strong and weak
analyticity). We now show that it takes
contractive values and that $T=M_S$.
Since ${\mathbf H}_2({\mathbb D})
\otimes {\mathcal H}$ is the reproducing kernel
Hilbert space with reproducing kernel
\[
\frac{I_{\mathcal H}}{1-\zeta\nu^*},
\]
convergence in norm implies
weak pointwise convergence in the coefficient space ${\mathcal
  H}$. Thus, for every $c\in{\mathcal H}$,
and every $f\in{\mathbf H}_2({\mathbb
D})\otimes{\mathcal H}$ in the form
\[
f(\zeta)=\sum_{n=0}^\infty \zeta^nf_n,
\]
using continuity of $T$, we have
\[
\begin{split}
\langle (Tf)(\zeta),c\rangle_{\mathcal H}&=
\langle (T(\sum_{n=0}^\infty M_\zeta^nf_n))(\zeta), c\rangle_{\mathcal H}
\\
&=\sum_{n=0}^\infty \langle (M_\zeta^nT(
f_n))(\zeta), c\rangle_{\mathcal H}
\\
&=\sum_{n=0}^\infty\langle \zeta^n (Tf_n)(\zeta),
c\rangle_{\mathcal H}\\
&=\sum_{n=0}^\infty
\langle \zeta^nS(\zeta)f_n,c\rangle_{\mathcal H}\\
&=\langle S(\zeta)f(\zeta)
c\rangle_{\mathcal H}.
\end{split}
\]
Thus $T=M_S$. Furthermore, the formula
\[
(M_S^*\frac{c}{1- \cdot\nu^*})(\zeta)=\frac{S(\nu)^*c}{1-\zeta\nu^*}
\]
implies that the kernel
\begin{equation}
\label{KS}
\frac{I_{\mathcal
H}-S(\zeta)S(\nu)^*}{1-\zeta\nu^*}
\end{equation}
is positive in ${\mathbb D}$ and in particular
$S$ takes contractive values.\\

STEP 2: {\sl The operators $M_{z_j}$ are bounded
in ${\mathcal H}(K_k)$.}\\

Indeed, let $g(z)=\sum_{\alpha\in\ell}c_\alpha
z^\alpha$ with
\[
\|g\|_{{\mathcal H}(K_k)}^2=
\sum_{\alpha\in\ell}|c_\alpha|^2 (2{\mathbb
N})^{-k\alpha}<\infty.
\]
Then, with $e_j\in\ell$ denoting the sequence
with all elements equal to $0$, at the exception
of the $j$-th, equal to $1$, we have
\[
z_jg(z)=\sum_{\alpha\in\ell}c_\alpha
z^{\alpha+e_j},\] and
\[
\begin{split}
\sum_{\alpha\in\ell} |c_\alpha|^2(2{\mathbb
N})^{-k(\alpha+e_j)} &=\sum_{\alpha\in\ell}
|c_\alpha|^2(2{\mathbb N})^{-k\alpha}
\left((2{\mathbb N})^{k\alpha}
(2{\mathbb N})^{-k(\alpha+e_j)}\right)\\
& \le\sum_{\alpha\in\ell} |c_\alpha|^2(2{\mathbb
N})^{-k\alpha},
\end{split}
\]
and so $M_{z_j}$ is a contractive operator in
${\mathcal H}(K_k)$.\\

STEP 3: {\sl We now take ${\mathcal H}
={\mathcal H}(K_k)$. The fact that $T$ commutes
with the operators of multiplication by $z_j$
implies that
\[
(T_n f)(z)=s_n(z)f(z),\]
for some ${\mathbf L}({\mathcal H}(K_k))$-valued
function $s_n(z)$.}\\

Indeed, we have in particular
\[
T_n(z^\alpha)=z^\alpha T_n(1).
\]
Thus, for $f(z)=\sum_{\alpha\in\ell}c_\alpha
z^\alpha\in{\mathcal H}(K_k)$ we have\\
\[
(T_nf)(z)= \sum_{\alpha\in\ell}c_\alpha T_n(1)=
T_n(1)f(z),
\]
first for finite sums, and then using continuity
of the operator $T_n$ for all $f\in{\mathcal
H}(K_k)$ in the norm topology, and finally
pointwise, since ${\mathcal H}(K_k)$ is a
reproducing kernel Hilbert space.\\

We set $(T_n(1))(z)=s_n(z)$. The result follows
with
\[
{\mathscr S}(\zeta,z)=\sum_{n=0}^\infty
\zeta^ns_n(z).
\]
Therefore $T$ is a multiplication operator in the
reproducing kernel Hilbert space ${\mathbf H}_2
({\mathbb D})\otimes{\mathcal H}(K_k)$. Since it
is a contraction,
\eqref{chaud_les_marrons} is in force.\\

The converse is clear: if a function $\mathscr
S$ is such that the kernel
\eqref{chaud_les_marrons} is positive on
${\mathbb
  D}\times{\mathbb K}_k$, then the
  operator of multiplication by
$\mathscr S$ is a contraction from ${\mathbf
H}_2({\mathbb D})\otimes{\mathcal H}(K_k)$ into
itself, and satisfies the commutativity
hypothesis \eqref{luxembourg}.
\mbox{}\qed\mbox{}\\

When the system is not random, ${\mathscr S}(\zeta,z)$
depends only on $\zeta$ and is a Schur
function.\\

We further remark that the function $S(\zeta)$ defined in \eqref{ST}
is a $ {\mathbf L}({\mathcal H}(K_k))$-valued
function, and the theory of realization for such
functions is thus applicable to it; see
\cite{adrs}. Using the results of \cite{adrs} we
have:

\begin{Tm}
\label{ma} Let $\mathscr S(\zeta, z)$ be a Schur
multiplier of the space ${\mathbf H}_2({\mathbb
D})\otimes{\mathcal H}(K_k)$, and let ${\mathcal
H}(\mathscr S)$ the associated reproducing
kernel Hilbert space. Then
\[
\mathscr S(\zeta,z)=D+\zeta C( I_{{\mathcal
H}(\mathscr S)}-\zeta A)^{-1}B,
\]
where
\[
\begin{pmatrix}A&B\\ C&D\end{pmatrix}\,\,:
\,\,\begin{pmatrix}{\mathcal
    H}(
\mathscr S)\\ {\mathcal H}(K_k)
\end{pmatrix}\longrightarrow
\begin{pmatrix}{\mathcal
    H}(
\mathscr S)\\ {\mathcal H}(K_k)\end{pmatrix}
\]
with
\[
\begin{split}
(Af)(\zeta, z)&=
\frac{f(\zeta, z)-f(0,z)}{\zeta},\\
(Bx)(\zeta, z)&=\frac{\mathscr
S(\zeta, z)-{\mathscr S}(0,z)}{\zeta}x(z),\\
(Cf)(\zeta, z)&=f(0,z),\\
(Dx)(\zeta, z)&=\mathscr S(0,z)x(z),
\end{split}
\]
defining a coisometric realization of $\mathscr S$.
\end{Tm}

{\bf Proof:} Indeed, let ${\mathcal H}$ denote a
Hilbert space, and let $S$ be an ${\mathbf
L}({\mathcal H})$-valued function analytic and contgractive in
${\mathbb D}$. Then the kernel \eqref{KS} is positive in the open unit
disk. Let ${\mathcal
  H}(S)$ be the associated reproducing kernel Hilbert space. The formulas

\[
\begin{split}
(Af)(\zeta)&=\frac{f(\zeta)-f(0)}{\zeta},\\
(B\xi)(\zeta)&=\frac{S(\zeta)-S(0)}{\zeta}\xi,\\
Cf&=f(0),\\
D\xi&=\mathscr S(0)\xi,
\end{split}
\]
define a coisometric realization with
\[
\begin{pmatrix}A&B\\ C&D\end{pmatrix}\,\,:\,\,\begin{pmatrix}{\mathcal
    H}(
 S)\\ {\mathcal H}\end{pmatrix}\longrightarrow
\begin{pmatrix}{\mathcal
    H}(
 S)\\ {\mathcal H}\end{pmatrix}
\]
in the form:
\[
S(\zeta)=D+\zeta C(I_{{\mathcal H}(S)}-\zeta
A)^{-1}B.
\]
Now, there is here an extra structure,
which does not appear in the general case: there is a function
$\mathscr S(\zeta, z)$ such that
\[
(S(\lambda)x)(z)=\mathscr S(\zeta, z)x(z),\quad x\in{\mathcal
  H}(K_k).
\]
The formulas for the realization of $\mathscr S$ then easily
follow.
\mbox{}\qed\mbox{}\\

Since constants belong to ${\mathbf H}_2
({\mathbb D})\otimes {\mathcal H}(K_k)$, it
follows in particular that ${\mathscr
S}(\zeta,z)\in {\mathbf H}_2({\mathbb D})\otimes
{\mathcal H}(K_k)$, and thus can be written as
\[
{\mathscr S}(\zeta,z)=
\sum_{n=0}^\infty\zeta^nh_n(z).
\]
The fact that the operator is a contraction is
then equivalent to the fact that the block
Toeplitz operator
\begin{equation}
\label{toep}
\begin{pmatrix}T_{h_0}&  0    &\cdot&\cdot\\
               T_{h_1}&T_{h_0}&\cdot&\cdot\\
               T_{h_2}&T_{h_1}&T_{h_0}&\cdot
                      &        & &
\end{pmatrix}
\end{equation}
is a contraction from the Hilbert space
$\ell_2({\mathbb N}_0)\otimes{\mathcal H}(K_k)$
into itself.\\

As a consequence, we have:

\begin{Tm}
Let $k>\l+1$ and let $(h_n)$ be a sequence of
elements in ${\mathcal H}_l$ indexed by ${\mathbb N}_0$. Let
\[
y_n=\sum_{m=0}^n h_{n-m}\lozenge u_m,\quad n=0,1,
\ldots
\]
Then there exists an $M\in[0,1]$ such that, for
all signals
\begin{equation}
\label{toeplitz}
\sum_{n=0}^\infty\|y_n\|^2_{k}\le M
\sum_{n=0}^\infty \|u_n\|^2_{k}
\end{equation}
if and only if \eqref{toep} is a contraction.
\end{Tm}

\section{The Bochner-Chandrasekharan theorem}
\setcounter{equation}{0}
\label{rotterdam}
Let ${\mathcal H}$ be a Hilbert space. We denote
for $h\in{\mathbb R}$ by $T_h$ the translation
operator:
\[
(T_hf)(t)=f(t+h),
\]
where $f$ is in ${\mathbf L}_2 ({\mathbb
R})\otimes{\mathcal H}$. The fact that a bounded
operator from ${\mathbf L}_2({\mathbb
R})\otimes{\mathcal H}$ into itself commutes
with all the $T_h$, expresses the time-invariance
of the underlying linear system. When ${\mathcal
H}={\mathbb C}$, the Bochner-Chandrasekharan
theorem provides a characterization of all bounded
operators from ${\mathbf L}_2({\mathbb R})$ into
itself which commute with the $T_h$; see
\cite[Theorem 72,
  p. 144]{boch_chan}. In the present
section we show how the strategy of
\cite{boch_chan} allows to prove a version of
this theorem in the case of Hilbert space valued
functions. Our approach consists of reducing the
operator-valued case to the scalar case, but we
need an extra assumption (the invariance property
\eqref{280909}), which means that we consider
causal systems, a property which holds in the
cases considered in the paper. We first prove the
following result. Its scalar version is
\cite[Theorem 70, p. 140]{boch_chan}. In the
statement of Theorem \ref{overwork}, for a given
Hilbert space
 ${\mathcal H}$, the Banach
space of bounded operator from ${\mathcal H}$
into itself is denoted by ${\mathbf L}({\mathcal
H})$. Moreover, we denote by ${\mathbf
H}_2({\mathbb C}_+)$ the Hardy space of the open
upper half plane ${\mathbb C}_+$.\\

In connection with this theorem, it is of
interest to mention the following result proved
in \cite[Theorem 1 p. 397]{MR0205028}: First
recall that ${\mathcal E}^\prime$ denotes the
space of distributions with compact support, and
that ${\mathcal D}^\prime$ denotes the dual of
the space ${\mathcal D}$ of infinitely
differentiable functions with compact support.
Then, any linear continuous map from ${\mathcal
E}^\prime$ into ${\mathcal D}^\prime$ which
commutes with the translation operators is a
convolution by a tempered distribution. See also
the discussion in \cite[p. 35]{lacroix-sonrier}.

\begin{Tm}
\label{overwork}
Let ${\mathcal H}$ be a separable Hilbert space,
and let $X$ be a bounded linear operator from
${\mathbf L}_2({\mathbb R})\otimes{\mathcal H}$
into itself. Assume that if $f$ vanishes outside
an interval, then $Xf$ vanishes outside the same
interval. Assume moreover that
\begin{equation}
X({\mathbf H}_2({\mathbb C}_+)\otimes {\mathcal
H}) \subset {\mathbf H}_2({\mathbb C}_+)\otimes
{\mathcal H}. \label{280909}
\end{equation}
Then, there exists a ${\mathbf L}({\mathcal
H})$-valued function $S$ bounded in norm and
analytic in ${\mathbb C}_+$ such that $X$ is the
operator of multiplication by $S$:
\[
Xf(\la)=S(\la)f(\la),\quad f\in{\mathbf
H}_2({\mathbb C}_+)\otimes{\mathcal
H},\quad\la\in {\mathbb C}_+.
\]
\end{Tm}

\begin{Rk}
Bochner and Chandrasekharan call an operator
which has the vanishing property of Theorem
\ref{overwork} a {\sl transformation of simple
type}.
\end{Rk}

{\bf Proof of Theorem \ref{overwork}:} Let
$e_0,e_1,\ldots$ be an orthonormal basis of
${\mathcal H}$. For $n,m\in{\mathbb N}_0$, the
expressions
\begin{equation}
\label{groningen_300908} (X_{n,m}(u))(t)=\langle
(X(ue_n))(t),e_m\rangle_{\mathcal H},\quad
u\in{\mathbf L}_2 ({\mathbb R}).
\end{equation}
define bounded linear operators from ${\mathbf
L}_2({\mathbb R})$ into itself since, by the
Cauchy-Schwarz inequality
\begin{equation*}
|(X_{n,m}(u))(t)|\le\|(X(ue_n))(t)\|_{\mathcal
H},
\end{equation*}
and so
\[
\begin{split}
\int_{\mathbb R} |(X_{n,m}(u))(t)|^2dt
&\le\int_{\mathbb R}\|X(ue_n)\|^2_{\mathcal H}
(t)dt\\
&=\|X(ue_n)\|^2_{{\mathbf L}_2({\mathbb
R})\otimes{\mathcal H}}
\\
&\le \|X\|^2\|ue_n\|^2_{{\mathbf L}_2({\mathbb
R})\otimes{\mathcal H}}\\
&=  \|X\|^2\cdot\|u\|^2_{{\mathbf L}_2({\mathbb
R})}.
\end{split}
\]
Moreover, the $X_{n,m}$ are transformations of
simple type since $X$ is a transformation of
simple type; therefore, by \cite[Theorem 70, p.
140]{boch_chan}, there exist bounded measurable
functions $s_{m,n}$ such that
\begin{equation}
\label{250909} \langle (X(ue_n))(t),e_m\rangle_{
{\mathcal H}}= s_{m,n}(t)\cdot u(t)\quad
u\in{\mathbf L}_2 ({\mathbb R}).
\end{equation}
From the proof in \cite{boch_chan} we moreover
have that
\[
\sup_{t\in{\mathbb R}}|s_{m,n}(t)|\le \|T\|.
\]
Setting $u=\langle f,e_n\rangle_{\mathcal H} $
with $f\in {{\mathbf L}_2({\mathbb R})\otimes
{\mathcal H}}$
  in \eqref{250909}, we obtain
\begin{equation}
\label{paris_texas!}
\big\langle (X(\langle f,e_n\rangle_{{\mathcal
H}} e_n))(t),e_m\big\rangle_{{\mathcal H}}=
s_{n,m}(t)\cdot \langle
f(t),e_n\rangle_{{\mathcal H}}
\end{equation}
for every $f\in{\mathbf L}_2 ({\mathbb
R})\otimes{\mathcal H}$.\\

Summing \eqref{paris_texas!} with respect to $n$
we expect to obtain
\begin{equation}
\label{oberkampf_ligne9} \big\langle
(Xf)(t),e_m\big \rangle_{\mathcal H}=
\sum_{n=0}^\infty s_{m,n}(t)\cdot \langle
f(t),e_n\rangle_{\mathcal H},\quad
m=0,1,2,\ldots,
\end{equation}
to show that
\[
(Xf)(t)=S(t)f(t),
\]
where $S(t)$ is the operator with matrix
$(s_{m,n}(t))$. We cannot justify these
computations in the general case (see also the
remark following the proof of the theorem), and
at this stage of the proof,
hypothesis \eqref{280909} is applied.\\

In view of \eqref{280909}, the functions
$s_{m,n}$ are analytic in the open upper
half-plane. The space ${\mathbf H}_2({\mathbb
C}_+)$ is the reproducing kernel Hilbert space
with reproducing kernel
\[ \frac{I_{\mathcal H}
}{-i(\lambda-\mu^*)},\quad \lambda,\mu\in{\mathbb
C}_+.
\]
and therefore convergence in norm in this
space implies pointwise convergence. Thus it
follows from
\begin{equation}
\label{vierzon} Xf=\sum_{0}^\infty \langle
Xf,e_n\rangle_{\mathcal H} e_n,
\end{equation}
where the convergence is in ${\mathbf H}_2
({\mathbb C}_+)\otimes {\mathcal H}$, that for
every $m\in{\mathbb N}_0$ and $\lambda
\in{\mathbb C}_+$,
\[
\begin{split}
\langle Xf,\frac{e_n}{-i(\cdot-\lambda^*)}
\rangle_{{\mathbf L}_2\otimes{\mathcal H}} &=
\langle (Xf)(\lambda),e_m\rangle_{\mathcal H}\\
&=\sum_{m=0}^\infty\langle \langle
(Xf)(\la),e_n\rangle_{\mathcal
H}e_n,e_m\rangle_{\mathcal H}\\
&=\sum_{m=0}^\infty s_{m,n}(\la)\langle
f(\la),e_n\rangle_{\mathcal H}\\
&=\langle S(\la)f(\la),e_n\rangle_{\mathcal H}.
\end{split}
\]
\mbox{}\qed\mbox{}\\

\begin{Rk}{\rm
Given a sequence of ${\mathbf L}_2({\mathbb R})$
functions which converges in ${\mathbf
L}_2({\mathbb R})$, there exists a subsequence
of this sequence which converges almost
everywhere; see for instance \cite[Th\'eor\`eme
2.3 p. 95]{descombes}. Therefore, without
assumption \eqref{280909}, convergence in
norm of the series \eqref{vierzon} and a
diagonal argument allows to assert that, for a
given $f\in{\mathbf L}_2({\mathbb R})
\otimes{\mathcal H}$, \eqref{oberkampf_ligne9}
holds almost everywhere on the real line when
the limit is taken via a subsequence.}
\end{Rk}

Recall that we denote by $U$ the Fourier
transform. We can now prove:
\begin{Tm}
Let $X$ be a bounded operator from ${\mathbf
L}_2({\mathbb R})\otimes{\mathcal H}$ into itself
which commute with translation operators, and
assume that $UXU^{-1}$ satisfies \eqref{280909}.
 Then there exists a ${\mathbf L}({\mathcal
H})$-valued function $S$ such that
\[
UXU^{-1}f=S(t)f(t), \quad f\in{\mathbf
H}_2({\mathbb C}_+)\otimes{\mathcal H}.
\]
\end{Tm}

{\bf Proof:} Let $X_{m,n}$ be defined by
\eqref{groningen_300908} and $h\in{\mathbb R}$.
Since $X$ commutes with the $T_h$, we have for
$u\in{\mathbf L}_2({\mathbb R}$
\[
\begin{split}
(T_h(X_{m,n}u))(t)&=\langle (X(ue_n))(t+h),
e_m\rangle_{\mathcal H}\\
&=\langle (X(T_hu))(t)e_m\rangle_{\mathcal H}\\
&=(X_{m,n}(T_hu))(t),
\end{split}
\]
and thus $X_{m,n}$ also commutes with the
translation operators. We now apply
\cite[Theorem 72, p. 144-147]{boch_chan} to show
that the image under the Fourier transform of
$X_{m,n}$ is of simple type, which allows to
conclude the thanks to the preceding Theorem \ref{overwork}.
\mbox{}\qed\\
\section{BIBO stability: The continuous time case}
\setcounter{equation}{0}
In the classical case, one considers the Hardy
space ${\mathbf H}_2({\mathbb C}_+)$ of the open
upper half-plane,
%
and, in the stochastic case, we will consider the
space ${\mathbf H}_2({\mathbb C}_+)\otimes
{\mathcal
H}_k$.\\

In this section we restrict ourselves to the
case of continuous functions $t\mapsto f(t)$
from the real line to ${\mathcal H}_k$ for some
$k\ge 1$. Then, there is no difficulty to define
the integral $\int_a^b f(t)dt$ on a compact
interval, and, with appropriate hypothesis, also
on the real line. Obviously, more general
situations can be considered. These will not be
pursued in
the present paper.\\

A continuous signal will be an ${\mathcal
H}_k$-valued continuous function $t\mapsto u(t)$
defined for $t\in{\mathbb R}$, and such that
\[
\sup_{t\in{\mathbb R}}\|u(t)\|_{k}<\infty.
\]
\begin{La}
Let $k>l+1$. Assume that $t\mapsto f(t)$ is a
continuous ${\mathcal H}_{l}$-valued function and
$t\mapsto g(t)$ is a continuous ${\mathcal
H}_{k}$-valued function. The ${\mathcal
H}_{k}$-valued function $t\mapsto f(t)\lozenge
g(t)$ is continuous with respect to the norm of
${\mathcal H}_k$.
\end{La}

{\bf Proof:} to prove the claim, it suffices to
write
\[
f(t_1)\lozenge g(t_1)-f(t_2)\lozenge g(t_2)=
(f(t_1)-f(t_2))\lozenge g(t_1)+f(t_2)\lozenge
(g(t_1)-g(t_2))
\]
for $t_1,t_2\in {\mathbb R}$ and use inequality
\eqref{vage}.
\mbox{}\qed\mbox{}\\

\begin{Dn}
A continuous linear time-invariant stochastic
system will be defined by its $S_{-1}$-valued
impulse response
\begin{equation}
\label{F} h(\tau,
\w)=\sum_{\alpha\in\ell}h_\alpha(\tau)
H_\alpha(\w).
\end{equation}
The associated input-output relation is in the
Wick convolution form
\[
y(t)=\int_{\mathbb R} h(t-\tau)\lozenge
u(\tau)d\tau
\]
\end{Dn}

\begin{Tm}
Let $k>l+1$. Let $t\mapsto h(t)$ be a continuous
${\mathcal H}_l$-valued function
$(t\in{\mathbb R})$. Then:\\
There  exists $M>0$ such that\\
$(a)$ The integrals
\begin{equation}
y(t)=\int_{\mathbb R}h(s)\lozenge
u(t-s)ds,\quad\forall t\in{\mathbb R},
\label{conv1}
\end{equation}
exist for all continuous ${\mathcal H}_k$-valued
functions $u$ such that
\[\sup_{\mathbb
R}\|u(t)\|_{k}<\infty\]
and\\
$(b)$ it holds that
\begin{equation}
\sup_{t\in\mathbb R}\|\int_{\mathbb
R}h(s)\lozenge u(t-s)ds\|_{k} \le
M\sup_{t\in\mathbb R}\|u(t)\|_{k}
\label{le_louvre}
\end{equation}
if and only if for every continuous $g\in
{\mathcal H}_k$ of norm less or equal to $1$, it
holds that
\begin{equation}
\label{nice}
\int_{\mathbb
R}\|T_{h(t)}^*g\|_{k}dt\le M
\end{equation}
\end{Tm}

{\bf Proof:} We first remark that the continuity
of the function $t\mapsto h(t)$ implies the
continuity of the function $t\mapsto
\|T^*_{h(t)}g\|_k$. Indeed, we have for
$t_1,t_2\in{\mathbb R}$
\[
\begin{split}
\big|\|T^*_{h(t_1)}g\|_k-\|T^*_{h(t_2)}g\|_k
\big|&\le
\|(T^*_{h(t_1)}-T^*_{h(t_2)})g\|_k\\
&\le
\|g\|_k\cdot\|T^*_{h(t_1)}-T^*_{h(t_2)}\|\\
&=\|g\|_k\cdot\|T^*_{h(t_1)-h(t_2)}\|\\
&=\|g\|_k\cdot\|T_{h(t_1)-h(t_2)}\|\\
&\le A(k-l)\|g\|_k\cdot\|{h(t_1)-h(t_2)}\|_l
\end{split}
\]
where we have used \eqref{belmondo} to obtain to the
last inequality. Thus, the integral on the left
handside of \eqref{nice} makes sense as a
possibly divergent integral. We now show that
it converges.\\

Assume that \eqref{le_louvre} is in force. Then,
for every real $t$,
\[
\int_{\mathbb R}\|h(\tau)\lozenge u(t-\tau)d\tau\|_{k}
\le M\sup_{t\in\mathbb R}\|u(t)\|_{k}.
\]
Hence, we have that
\[
\int_{\mathbb R} |\langle h(\tau)\lozenge u(t-\tau)d\tau,
g\rangle_{k}|du\le M \sup_{t\in\mathbb
R}\|u(t)\|_{k},
\]
and thus
\[
|\int_{\mathbb R} \langle h(s)\lozenge u(t-s),
g\rangle_{k}ds|\le M \sup_{\mathbb
R}\|u(t)\|_{k},
\]
for every $g\in {\mathcal H}_k$ of norm less or
equal to $1$. Hence, for such $g$,
\begin{equation}
\label{yes} |\int_{\mathbb R} \langle u(t-s),
T_{h(s)}^*(g)\rangle_{k}ds|\le M \sup_{\mathbb
R}\|u(t)\|_{k}.
\end{equation}
Let $\epsilon>0$ and $t$ fixed in ${\mathbb R}$.
The function
\[
u(s)=\frac{T_{h(t-s)}^*(g)}{\|T_{h(t-s)}^*(g)
\|_{k}+\epsilon}
\]
is continuous in $s$ with respect to the norm of
${\mathcal H}_k$. Its substitution in
\eqref{yes} leads to
\[
\int_{\mathbb R}\frac{\|T^*_{h(s)}(g)\|_k^2}{
\|T^*_{h(s)}(g)\|_k+\epsilon}ds\le M.
\]
Taking $\epsilon=1/n$, $n=1,2,\ldots$ and using
the Monotone Convergence Theorem then
leads to \eqref{nice}.\\

Conversely, assume that \eqref{nice} holds, and
let $u$ be a ${\mathcal H}_k$-valued continuous
and bounded function. From
\[
\begin{split}
|\langle h(s)\lozenge u(t-s), g\rangle_{k}|&=
|\langle u(t-s), T_{h(s)}^*(g)\rangle_{k}|\\&\le
\sup_{s\in\mathbb
R}\|u(s)\|_{k}\|T_{h(s)}^*(g)\|_{k},
\end{split}
\]
we obtain that
\[
\int_{\mathbb R}|\langle h(s)\lozenge u(t-s),
g\rangle_{k}|ds\le M\sup_{s\in\mathbb
R}\|u(s)\|_{k},
\]
and thus
\[
\int_{\mathbb R}|\langle h(s)\lozenge u(t-s),
g\rangle_{k}|ds\le M\sup_{s\in\mathbb
R}\|u(s)\|_{k}.
\]
Since this inequality holds for all $g\in
{\mathcal H}_k$ of norm less or equal to $1$, we
obtain \eqref{le_louvre}. In the above chain of
equalities we have used the following (see
\cite[(8.7.6) p. 169]{dieudonne3}): given a
Hilbert space ${\mathcal H}$ and a ${\mathcal
H}$-valued continuous function $h$ such that
$\int_{\mathbb R} h(s)ds$ exists, then,  for any
$g\in{\mathcal H}$, it holds that
\[
\langle \int_{\mathbb R} h(s)ds, g
\rangle_{\mathcal H}=\int_{\mathbb R}\langle
h(s),g\rangle_{\mathcal H}ds,
\]
which concludes the proof of the theorem.
\mbox{}\qed\mbox{}\\

We note that a weaker condition than
\eqref{nice} is given by

\begin{equation}
\label{nice1} \int_{\mathbb R}\|T_{h(t)}\|dt\le
M,
\end{equation}
where, in view of the continuity of $h$, the
integral makes sense since, by \eqref{belmondo},
\[
\|T_{h(t_1)-h(t_2)}\| \le
A(k-l)\|{h(t_1)-h(t_2)}\|_l.
\]
\mbox{}\\

We note that, when $h$ is not random, the Wick
product reduces to a pointwise product, and so
\eqref{nice} reduces to the well known condition
\[
\int_{\mathbb R}\|h\|(t)dt<\infty.\]
See for instance \cite[\S2.6.1. p. 175]{MR569473}
or \cite[Corollary 3, p. 585]{MR0349260}.\\

Finally, write $u(t)=\sum_\alpha
u_\alpha(t)H_\alpha$, and similarly for the
output. Taking the Laplace transform (denoted by
$\widehat{u}$), followed by the Hermite
transform, we obtain an expression of the form
\[
\sum_{\alpha\in\ell} \widehat{u_\alpha}(\xi)
z^\alpha,
\]
where we have denoted by $\xi$ the variable
corresponding to the Laplace transform.

\begin{Pn} It holds that
\begin{equation}
\label{Richard-Lenoir} \widehat{
y}(z,\zeta)=\sum_{\alpha\in\ell}z^\alpha
  \left\{\sum_{\beta\le\alpha}
\widehat{h}_\beta(\zeta){\widehat
u}_{\alpha-\beta}(\zeta)\right\}.
\end{equation}
\end{Pn}
{\bf Proof:} For given real numbers $\tau$ and
$t$ we have
\[
h(t-\tau)\lozenge u(\tau)=\sum_{\alpha\in\ell}
\left\{\sum_{\beta\le
  \alpha}
h_\beta(t-\tau)u_{\alpha-\beta}(\tau)\right\}
H_\alpha(\w),
\]
and so, $y_\alpha$ is given by
\begin{equation}
\label{Breguet-Sabin}
y_\alpha(\tau)=\int_{\mathbb R} \sum_{\beta\le
  \alpha}h_\beta(t-\tau)u_{\alpha-\beta}(\tau)
  d\tau
\end{equation}
Taking Hermite and Laplace transforms leads to
\eqref{Richard-Lenoir}. \mbox{}\qed\mbox{}\\

We see that as in the discrete time case, in
continuous time, we also have two convolutions,
one with respect to time and one reflecting the
stochastic aspect of the system.

\section{Dissipative continuous time
random systems} \setcounter{equation}{0} For
$t\ge 0$ let $\chi_t$ denote the function
\[
\chi_t(\la)= e^{it\lambda}.
\]
Note that, for $t\ge 0$, the operator
$M_{\chi_t}$ of multiplication by $\chi_t$ is an
isometry from ${\mathbf H}_2({\mathbb C}_+)$
into itself, and that, by the
Bochner-Chandrasekharan theorem and the
properties of the Fourier transform, a bounded
linear operator from ${\mathbf H}_2({\mathbb
C}_+)$ into itself commutes with the operators
$M_{\chi_t}$ if and only if it is a
 multiplication
operator, or, in systems terminology, if and only
if it is time-invariant.

\begin{Tm}
Let $T$ be a linear contractive operator from
${\mathbf H}_2({\mathbb C}_+)\otimes {\mathcal
H}(K_k)$ into itself and such that
\begin{equation}
\label{230909}
\begin{split}
T(\chi_tf)&=\chi_tf\\
T(z_jf)&=z_jTf.
\end{split}
\end{equation}
Then, $T$ is a multiplication operator by a
function ${\mathscr H} (\lambda, z)$ which is
such that the kernel
\begin{equation}
\label{vezoul} (1-{\mathscr
H}(\lambda,z){\mathscr H}(\nu,w)^*)
\frac{K_k(z,w)}{-i(\lambda-\nu^*)}
\end{equation}
is positive in ${\mathbb C}_+\times {\mathbb
K}_k$.
\end{Tm}

{\bf Proof:} The proof follows the lines of the
proof of Theorem \ref{KSU_080909}. We use the
Hilbert space version of the
Bochner-Chandrasekharan theorem to assert the
existence of a ${\mathbf L}({\mathcal
H}(K_k))$-valued function $S$ such that
\[
T(\lambda)=S(\lambda)f(\lambda),\quad f\in
{\mathbf H}_2({\mathbb C}_+)\otimes{\mathcal
H}(K_k),
\]
for which  $\|S(\lambda)\|_{{\mathbf
L}({\mathcal H}(K_k))} \le 1,\quad \forall
\lambda\in{\mathbb C}_+$. This last norm
condition on $S$ is equivalent to the positivity
of the kernel
\[
\frac{I_{{\mathcal H}(K_k)}-
S(\lambda)S(\nu)^*}{-i(\lambda- \nu^*)}
\]
in the open upper half-plane ${\mathbb C}_+$. Let
$\lambda\in {\mathbb C}_+$.
%
Since the operator
\[
S(\la)\,:\,\, {\mathcal H}(K_k)\longrightarrow
 {\mathcal H}(K_k)
 \]
commutes with the operators of multiplication
 by the variables $z_j$ we have that, for every
 $\alpha\in\ell$,
\begin{equation}
\label{ribambelle}
((S(\la))(M_{z^{\alpha}}))(z)=z^\alpha
((S(\la)(1)))(z).
\end{equation}
We set:
\[
{\mathscr S}(\lambda, z)=S(\la)(1).
\]
Let $f\in{\mathcal H}(K_k)$. It follows from
\eqref{ribambelle} that
\[
(S(\lambda)(f))(\la, z)={\mathscr S}(\la, z)\cdot
f(\la, z). \]

The operator $T$ is thus a contractive
multiplication operator in the reproducing
kernel Hilbert space ${\mathbf H}_2({\mathbb
C}_+) \otimes {\mathcal H}(K_k)$, and the kernel
\eqref{vezoul} is positive in ${\mathbb C}_+\times {\mathbb K}_k$.
\mbox{}\qed\mbox{}\\

As in the disk case, one can write a realization
for $\mathscr S$ in terms of an associated
reproducing kernel Hilbert space. This will not
be pursued here.\\

When the system is not random, setting $z=w$ and
$\zeta=\nu$ in \eqref{vezoul}, we get that the
function ${\mathscr S}(\zeta)$ is
contractive in ${\mathbb C}_+$.\\

\section{${\mathbf L}_2$-${\mathbf L}_\infty$
stability} In this section we prove the analog
of Theorem \ref{salomon} in the continuous time
case.

\begin{Dn} Let $k>l+1$.
The  continuous time system \eqref{conv1} (with
a continuous $h$) is said to be ${\mathbf
L}_2$-${\mathbf L}_\infty$ stable if there
exists an $M<\infty$ such that
\begin{equation}
\sup_{t\in{\mathbb R}}\|\int_{\mathbb R}h(s)\lozenge
u(t-s)ds\|_k\le M
\left(\int_{\mathbb
  R}\|u(t)\|^2dt\right)^{1/2}
\label{jaune}
\end{equation}
for all ${\mathcal H}_k$-valued functions $u$ which are
moreover continuous in norm and for which the right handside of
\eqref{jaune}
is finite.
\end{Dn}

\begin{Tm}
The system \eqref{conv1} is
${\mathbf L}_2$-${\mathbf L}_\infty$ stable if and only if
\begin{equation}
\label{opera_garnier}
\sup_{\substack{g\in{\mathcal H}_k\\
\|g\|_k\le 1}}\int_{\mathbb R}\|T_{h(t)}^*g\|_k^2dt<\infty.
\end{equation}
\end{Tm}
{\bf Proof:} We first assume that \eqref{opera_garnier} is in
force. Let $g\in{\mathcal H}_k$ of norm less or equal to $1$.
Using the Cauchy-Schwarz inequality, we obtain
\[
\begin{split}
|\int_{\mathbb R}\langle h(s)\lozenge u(t-s),g\rangle_{{\mathcal
 H}_k}ds|&=
|\int_{\mathbb R}\langle u(t-s),T_{h(s)}^*g\rangle_{{\mathcal H}_k}ds|\\
&\le\int_{\mathbb R}
\|T_{h(s)}^*g\|_k\|u(t-s)\|_kds\\
&\le \left(\int_{\mathbb R}\|T_{h(s)}^*g\|_k^2ds\right)^{1/2}\left(
\int_{\mathbb R}\|u(t-s)\|_k^2ds\right)^{1/2}\\
&=\left(\int_{\mathbb R}\|T_{h(s)}^*g\|_k^2ds\right)^{1/2}\left(
\int_{\mathbb R}\|u(s)\|_k^2ds\right)^{1/2},
\end{split}
\]
and thus the system is ${\mathbf L}_2$-${\mathbf L}_\infty$ stable
since
\[
\|\int_{\mathbb R}h(s)\lozenge
u(t-s)ds\|_k=
\sup_{\substack{g\in{\mathcal H}_k\\
\|g\|_k\le 1}}\int_{\mathbb R}
\langle h(s)\lozenge u(t-s),
g\rangle_{{\mathcal H}_k}
\]
Conversely, assume that the system \eqref{conv1}
is ${\mathbf L}_2$-${\mathbf L}_\infty$ stable.
Then for every $g$ as above,
\begin{equation}
\label{12345}
\begin{split}
|\int_{\mathbb R}\langle
h(s)\lozenge u(t-s),g\rangle_{{\mathcal
 H}_k}ds|&=
|\int_{\mathbb R}\langle u(t-s),
T_{h(s)}^*g\rangle_{{\mathcal
 H}_k}ds|\\
&\le M\left(\int_{\mathbb
R}\|u(t)\|_k^2dt\right)^{1/2}
\end{split}
\end{equation}
and so for every $t$, and in particular for $t=0$, the functional
\[
u\mapsto
\int_{\mathbb R}\langle u(t-s),
T_{h(s)}^*g\rangle_{{\mathcal H}_k}ds
\]
is bounded. It follows from Riesz
representation theorem for bounded
functionals on Hilbert space that
\[
\int_{\mathbb R}\|T_{h(s)}^*g\|_k^2ds<\infty.
\]
This inequality is moreover uniform in $g$ over all 
functions $g$ such
that $\|g\|_k\le 1$, since $g$ does not appear on
the right handside of the last inequality in
\eqref{12345}.
\mbox{}\qed\mbox{}\\

We note that condition \eqref{opera_garnier}
will hold as soon as
\[
\int_{\mathbb R}\|T_{h(s)}\|^2ds<\infty.
\]
In the nonrandom case, this means that $h$ is in
${\mathbf L}_2({\mathbb R})$,  and we obtain the
continuous time analog
of Theorem \ref{l1l2}. See
\cite[Table 2.2 p. 19 and pp. 21-25]{MR1200235}.\\

In connection with the preceding result, we
recall the following inequalities: First, if $p$
and $q$ are positive real numbers such that
$1/p+1/q=1$, and $f\in{\mathbf L}_p({\mathbb
R})$ and $g\in{\mathbf L}_q({\mathbb R})$, we
have that the convolution $f*g\in{\mathbf
L}_\infty({ \mathbb R})$ and
\[
\|f*g\|_\infty\le \|f\|_p\cdot \|g\|_q,\]
where $\|\cdot\|_p$ denotes the norm in ${\mathbf
L}_p({\mathbb R})$; see \cite[Th\'eor\`eme 2.3,
p. 148]{descombes}, and \cite[Corollary 1 p. 280]{MR37:726} for a more
general statement. The situation at hand here
corresponds to $p=q=2$. We also recall
 (see
\cite[Th\'eor\`eme IV.15, p. 66]{Brezis},
\cite[Exercice 4, p. 141]{Rudin-french}): If
$f\in {\mathbf L}_1({\mathbb R})$ and
$g\in{\mathbf L}_p({\mathbb R})$, with
$p\in[1,\infty]$ ($\infty$ included), then the
convolution $f*g \in{\mathbf L}_p({\mathbb R})$,
and
\[
\|f*g\|_p\le \|f\|_1\cdot\|g\|_p.
\]
These results should lead to other stability
results in the random case. Such
a line of research will not be developed here.\\

 {\bf Acknowledgments:} Daniel Alpay wishes to
thank the Earl Katz family for endowing the
chair  which supported his research. It is a
pleasure to thank Prof. L. Baratchart for
pointing out to us reference \cite{MR1200235}
and to Prof. Juliette Leblond for a discussion
related to the continuous time counterpart of
Theorem
\ref{l1l2}.\\

\bibliographystyle{plain}
\def\cprime{$'$} \def\lfhook#1{\setbox0=\hbox{#1}{\ooalign{\hidewidth
  \lower1.5ex\hbox{'}\hidewidth\crcr\unhbox0}}} \def\cprime{$'$}
  \def\cprime{$'$} \def\cprime{$'$} \def\cprime{$'$} \def\cprime{$'$}

\end{document}